\documentclass[12pt,a4paper]{article}
\usepackage[cp1251]{inputenc} 
\usepackage{amsfonts}

\newcommand{\di}{\displaystyle}

\newcommand{\al}{\alpha}
\newcommand{\be}{\beta}
\newcommand{\ga}{\gamma}
\newcommand{\de}{\delta}
\newcommand{\la}{\lambda}
\newcommand{\om}{\omega}

\newcommand{\vv}{\varphi}
\newcommand{\iy}{\infty}

\begin{document}

\begin{center}
{\large\bf
On Quasi-periodic Differential Pencils with Jump Conditions Inside the Interval}\\[0.1cm]
{\bf V.\,Yurko} \\[0.1cm]
\end{center}

\thispagestyle{empty}

{\bf Abstract.} Non-self-adjoint second-order differential pencils on a finite interval
with non-separated quasi-periodic boundary conditions and jump conditions are studied.
We establish properties of spectral characteristics and investigate the inverse spectral
problem of recovering the operator from its spectral data. For this inverse problem
we prove the corresponding uniqueness theorem and provide an algorithm for
constructing its solution.

\noindent {\bf 1. Introduction}\\

Consider the boundary value problem $B$ of the form
$$
y''+(\rho^2+\rho p(x)+q(x))y=0,\quad x\in [0,T],                                  \eqno(1)
$$
$$
y(0)=\al y(T),\quad y'(0)-(i\rho h'+h)y(0)=\be y'(T),                             \eqno(2)
$$
$$
y(b_j+0)=\ga_j y(b_j-0),\;
y'(b_j+0)=\ga_j^{-1}y'(b_j-0)+(i\rho\eta'_j+\eta_j)y(b_j-0),\;j=\overline{1,N-1}, \eqno(3)
$$
$$
0=b_0<b_1<\ldots <b_{N-1}<b_N=T,
$$
where $\rho$ is the spectral parameter, $p(x), q(x)$ are complex-valued
functions, and $h', h, \al, \be, \ga_j$, $\eta'_j, \eta_j$ are complex numbers,
$\al\be\ga_j\ne 0,$ $z_0^{\pm}:=\al(1\mp h')+\be\ne 0,$
$\xi_j^{\pm}:=(\ga_j+(\ga_j)^{-1})/2\mp\eta'_j/2\ne 0.$
Assume that $p(x)\in AC[0,T]$ and $q(x)\in L(0,T).$
We study a nonlinear inverse problem of recovering $B$ from its spectral data.
For this inverse problem we prove the uniqueness theorem and provide
a procedure for constructing its solution.

Inverse spectral problems often appear in mathematics as well as in applications [1-3].
For {\it Sturm-Liouville operators} with separated boundary conditions, inverse spectral
problems have been studied fairly completely (see the monographs [1-3] and the references
therein). Such problems for Sturm-Liouville operators with non-separated boundary conditions
investigated in [4-7] and other works.

Differential equations with nonlinear dependence on the spectral parameter arise
in various problems of mathematics as well as in applications. In particular,
several examples of such spectral problems arising in mechanical engineering are
provided in the book [8] of Collatz; see also [9-11], where further references
and links to applications can be found. Detailed studies on direct spectral
problems for some classes of ordinary differential operators depending nonlinearly
on the spectral parameter can be found in various publications, see e.g. [9-11].

{\it Inverse spectral problems} for differential pencils, because of their nonlinearity,
are more difficult for investigating, and nowadays there are only isolated
fragments, not constituting a general picture, in the inverse problem theory
for equation (1). Some aspects of the inverse problem theory for pencils under
various restrictions were studied in [12]-[16] and other works, but mostly
only particular questions are considered there. Inverse problems for general
non-selfadjoint boundary value problem (1)-(3) with jump conditions inside the interval
have not been studied yet. We note the inverse problem, considered in this paper,
appears in the inverse problem theory for differential operators on spatial networks
with cycles (see [17-19]) which have many applications in natural sciences and engineering.

Some words about the structure of the paper. The statement of the inverse problem
is provided in section 2. In Section 3 properties of the spectrum are established.
In particular, the Weyl-type function and the corresponding Weyl sequence are
introduced and investigated. In Section 4 we provide the solution of the inverse
spectral problem for the boundary value problem $B.$\\

{\bf 2. Statement of the inverse problem.}\\

Denote by $S(x,\rho)$ and $C(x,\rho)$ the solutions of equation (1) satisfying
jump conditions (3) and the initial conditions
$$
S(0,\rho)=C'(0,\rho)=0,\quad S'(0,\rho)=C(0,\rho)=1.
$$
For each fixed $x,$ the functions $S^{(\nu)}(x,\rho)$ and $C^{(\nu)}(x,\rho),\;
\nu=0,1,$ are entire in $\rho$ of exponential type, and
$\langle C(x,\rho), S(x,\rho)\rangle\equiv 1,$
where $\langle y, z\rangle:=yz'-y'z$ is the Wronskian of $y$ and $z.$
Put $\vv(x,\rho)=C(x,\rho)+(i\rho h'+h)S(x,\rho),$ $d(\rho)=S(T,\rho),$
$d_1(\rho)=C(T,\rho).$ Eigenvalues ${\cal P}=\{\rho_n\}_{n\in{\bf Z}}$ of
the boundary value problem (1)-(3) coincide with the zeros (counting with
multiplicities) of the characteristic function
$$
a(\rho)=\al\vv(T,\rho)+\be S'(T,\rho)-(1+\al\be).                            \eqno(4)
$$
Let $\Lambda:=\{n:\; n=\pm 1,\pm 2,\ldots\}={\bf Z}\setminus\{0\},$ and let
${\cal V}=\{\nu_n\}_{n\in\Lambda}$ be zeros (counting with multiplicities)
of $d(\rho).$ Then $\{\nu_n\}_{n\in\Lambda}$ are the eigenvalues of the
boundary value problem ${\cal B}$ for Eq. (1) with jump conditions (3)
and with the boundary conditions $y(0)=y(T)=0.$ The function $d(\rho)$
is called the characteristic function for ${\cal B}$.
Without loss of generality, we agree that the numeration is chosen such
that $\nu_n\ne\nu_k$ if $nk<0.$ Let $m_n$ be the multiplicity of $\nu_n$
($\nu_n=\nu_{n+1}=\ldots =\nu_{n+m_n-1}$). Put
$I:=\{n\in\Lambda:\;\nu_{n-1}\ne\nu_n\},$ $I':=\{n\in I:\;m_n>1\}.$

Denote $D(\rho)=\al\vv(T,\rho)+\be S'(T,\rho),\;
Q(\rho)=\al\vv(T,\rho)-\be S'(T,\rho).$ Then
$$
D(\rho)=d(\rho)+(1+\al\be),                                                \eqno(5)
$$
$$
\vv(T,\rho)=\frac{1}{2\al}\Big(D(\rho)+Q(\rho)\Big),\quad
S'(T,\rho)=\frac{1}{2\be}\Big(D(\rho)-Q(\rho)\Big).
$$
Since $\vv(x,\rho)S'(x,\rho)-\vv'(x,\rho)S(x,\rho)\equiv 1,$ it follows that
$$
Q^2(\rho)=D^2(\rho)-4\al\be(1+\vv'(T,\rho)S(T,\rho)),                      \eqno(6)
$$
and consequently,
$$
\dot Q(\rho)Q(\rho)=\dot D(\rho)D(\rho)-
2\al\be(\dot\vv'(T,\rho)S(T,\rho)+\vv'(T,\rho)\dot S(T,\rho)),             \eqno(7)
$$
where "dot" denotes derivatives with respect to $\rho.$
Let $n\in I.$ Denote
$$
\omega_n=\left\{ \begin{array}{ll}
0,\quad & Q(\nu_n)=0, \\
+1,\quad & Q(\nu_n)\ne 0,\; \mbox{arg}\, Q(\nu_n)\in[0,\pi), \\
-1,\quad & Q(\nu_n)\ne 0,\; \mbox{arg}\, Q(\nu_n)\in[\pi,2\pi),
\end{array}\right.
$$
$\omega_{n\nu}:=d_1^{(\nu)}(\nu_n),$ $\nu=\overline{0,m_n-1},$
$I_0=\{n\in I':\; \om_n=0\},$ $I_1=\{n\in I':\; \om_n\ne 0\}.$
The sequence $\Omega=\{\omega_n\}_{n\in I}\cup\{\omega_{n\nu}\}_{n\in I_0,
\;\nu=\overline{1,m_n-1}}$ is called the $\Omega$- sequence for $L.$
We note that if $I'=\emptyset$ (i.e. $m_n=1$ for all $n$), then
$\Omega=\{\omega_n\}_{n\in\Lambda}$.

Let $\al, \be$ and $\ga_j$ are known a priori and fixed.
The inverse problem is formulated as follows.

\smallskip
{\bf Inverse problem 1. }
Given $a(\rho), d(\rho)$ and $\Omega,$ construct $B.$

\smallskip
Obviously, in general it is not possible to recover all coefficients
from (2)-(3). Note that this inverse problem is a generalization of the
classical inverse problems for Sturm-Liouville operators [4-7].\\

{\bf 3. Properties of the spectral characteristics.}\\

Let $\Phi(x,\rho)$ be the solution of equation (1) under the jump
conditions (3) and the boundary conditions $\Phi(0,\rho)=1,\;\Phi(T,\rho)=0.$
Denote $M(\rho):=\Phi'(0,\rho).$ The function $M(\rho)$ is called the
Weyl-type function. Clearly,
$$
\Phi(x,\rho)=C(x,\rho)+M(\rho)S(x,\rho),                                      \eqno(8)
$$
$$
M(\rho)=-\frac{d_1(\rho)}{d(\rho)}.                                           \eqno(9)
$$
Since $\langle C(x,\rho), S(x,\rho)\rangle\equiv 1,$ it follows from (8) that
$$
\langle\Phi(x,\rho), S(x,\rho)\rangle\equiv 1.                               \eqno(10)
$$

Denote $T_k:=b_k-b_{k-1},\; k=\overline{1,N}.$ Then $b_k=T_1+\ldots+T_k$,
$T=T_1+\ldots+T_N$. Put $x_k=x-b_{k-1}$ for $x\in[b_{k-1}, b_k];$ hence
$x_k\in[0,T_k].$ Let $S_k(x_k,\rho)$ and $C_k(x_k,\rho)$ be the solutions
of equation (1) on $[b_{k-1}, b_k]$ under the initial conditions
$$
S_k(0,\rho)=C'_k(0,\rho)=0,\quad S'_k(0,\rho)=C_k(0,\rho)=1.                \eqno(11)
$$
For each fixed $x_k$, the functions $S^{(\nu)}_k(x_k,\rho)$ and
$C^{(\nu)}_k(x_k,\rho),\; \nu=0,1.$ are entire in $\rho,$ and
$$
\langle C_k(x_k,\rho), S_k(x_k,\rho)\rangle\equiv 1.
$$

{\bf Lemma 1. }{\it The following relations hold for
$k=\overline{1,N-1},\;\nu=0,1$:}
$$
S^{(\nu)}(b_{k+1}-0,\rho)=\ga_k S(b_k-0,\rho)C_{k+1}^{(\nu)}(T_{k+1},\rho)
+\ga_k^{-1}S'(b_k-0,\rho)S_{k+1}^{(\nu)}(T_{k+1},\rho)
$$
$$
+(i\rho\eta'_k+\eta_k) S(b_k-0,\rho)S_{k+1}^{(\nu)}(T_{k+1},\rho),         \eqno(12)
$$
$$
C^{(\nu)}(b_{k+1}-0,\rho)=\ga_k C(b_k-0,\rho)C_{k+1}^{(\nu)}(T_{k+1},\rho)
+\ga_k^{-1}C'(b_k-0,\rho)S_{k+1}^{(\nu)}(T_{k+1},\rho)
$$
$$
+(i\rho\eta'_k+\eta_k) C(b_k-0,\rho)S_{k+1}^{(\nu)}(T_{k+1},\rho),         \eqno(13)
$$

Indeed, fix $k=\overline{1,N-1}.$ Let $x\in[b_k,b_{k+1}],$ i.e.
$x=x_{k+1}+b_k$, $x_{k+1}\in[0,T_{k+1}].$ Using the fundamental
system of solutions $C_{k+1}(x_{k+1},\rho), S_{k+1}(x_{k+1},\rho),$ one has
$$
S^{(\nu)}(x,\rho)=A(\rho)C^{(\nu)}_{k+1}(x_{k+1},\rho)
+B(\rho)S^{(\nu)}_{k+1}(x_{k+1},\rho),\quad \nu=0,1.
$$
Taking initial conditions (11) into account we find the coefficients
$A(\rho)$ and $B(\rho),$ and arrive at (12). Relation (13) is proved similarly.

Denote
$$
{\cal E}(x)=\di\frac{1}{2}\di\int_0^x p(t)\,dt,\;
\om=\frac{1}{2T}\int_0^T p(t)\,dt,\; \tau=\mbox{Im}\,\rho,
\; G_{\de}=\{\rho:\; |\rho-\nu_n|\ge\de\;\forall n\},
$$
$$
\Pi^{\pm}=\{\rho:\; \pm\tau>0\},\;
\Pi^{+}_\de=\{\rho:\; \mbox{arg}\,\rho\in[\de,\pi-\de]\},\;
\Pi^{-}_\de=\{\rho:\; \mbox{arg}\,\rho\in[\pi+\de,2\pi-\de]\}.
$$
It is known (see [9]) that there exists a fundamental system of
solutions \\ $\{Y_1^{\pm}(x,\rho),\; Y_2^{\pm}(x,\rho)\},\; x\in[0,T],
\; \rho\in\Pi^{\pm},$ of equation (1) with the properties:\\
1) The functions $Y_k^{\pm}(x,\rho)$ are regular in
$\rho\in\Pi^{\pm},\; |\rho|>\rho^{*},$ and are continuous for
$x\in [0,T],\; \rho\in\overline{\Pi^{\pm}},\; |\rho|\ge\rho^{*}$.\\
2) For $|\rho|\to\iy,\; \rho\in\overline{\Pi^{\pm}},\; k=1,2,\; \nu=0,1,$
$$
\frac{d^\nu}{dx^\nu}Y_k^{\pm}(x,\rho)=(\rho R_k)^{\nu}
\exp(\rho x+{\cal E}(x))R_k)[1],\; [1]=1+O(\rho^{-1}),\; R_k=(-1)^{k-1}i.     \eqno(14)
$$

Using (14) and Lemma 1, one gets for $x\in(b_j,b_{j+1}),$
$|\rho|\to\iy,\;\rho\in\Pi_\de^{\pm}$:
$$
C^{(\nu)}(x,\rho)=
\frac{\xi_1^{\pm}\ldots \xi_j^{\pm}}{2}
(\mp i\rho)^{\nu}\exp(\mp i(\rho x+{\cal E}(x))[1],                          \eqno(15)
$$
$$
S^{(\nu)}(x,\rho)=\mp\frac{\xi_1^{\pm}\ldots \xi_j^{\pm}}{2i\rho}
(\mp i\rho)^{\nu}\exp(\mp i(\rho x+{\cal E}(x))[1],                          \eqno(16)
$$
$$
\Phi^{(\nu)}(x,\rho)=\frac{1}{\xi_1^{\pm}\ldots \xi_j^{\pm}}
(\pm i\rho)^{\nu}\exp(\pm i(\rho x+{\cal E}(x))[1],                          \eqno(17)
$$
In particular, we have for $|\rho|\to\iy,\;\rho\in\Pi_\de^{\pm}$:
$$
a(\rho)=\frac{z_0^{\pm}}{2}(\xi_1^{\pm}\ldots \xi_{N-1}^{\pm})
\exp(\mp i(\rho+\om)T)[1],\;
d(\rho)=\mp\frac{1}{2i\rho}(\xi_1^{\pm}\ldots \xi_{N-1}^{\pm})
\exp(\mp i(\rho+\om)T)[1].                                                  \eqno(18)
$$
Moreover, for $x\in[0,T],\; \rho\in\overline{\Pi^{\pm}}$:
$$
|S^{(\nu)}(x,\rho)|\le C|\rho|^{\nu-1}\exp(|\tau|x),\;
|C^{(\nu)}(x,\rho)|\le C|\rho|^{\nu}\exp(|\tau|x),                          \eqno(19)
$$
$$
|\Phi^{(\nu)}(x,\rho)|\le C |\rho|^{\nu}\exp(-|\tau|x),
\; |M(\rho)|\le C|\rho|,\; \rho\in G_{\de}.                                 \eqno(20)
$$

Let $n\in I.$ Using (9) we obtain that in a neighborhood of the
point $\rho=\nu_n$, the function $M(\rho)$ has the representation
$$
M(\rho)=\sum_{\nu=0}^{m_n-1}
\frac{M_{n+\nu}}{(\rho-\nu_n)^{\nu+1}}+M^*_n(\rho),                         \eqno(21)
$$
where $M^*_n(\rho)$ is regular in $\rho=\nu_n$, and the coefficients
$M_{n+\nu},\; \nu=\overline{0,m_n-1}$ are calculated from
$d_1^{(\nu)}(\nu_n)$ and $d^{(\nu+m_n)}(\nu_n)$ for
$\nu=\overline{0,m_n-1}.$ More precisely,
$$
M_{n+m_n-1-\nu}=-\frac{1}{d_{0n}}\Big(d^1_{\nu n}+\sum_{k=0}^{\nu-1}
M_{n+m_n-1-k}d_{\nu-k,n}\Big),\quad \nu=\overline{0,m_n-1},
$$
$$
d^1_{\nu n}:=\frac{1}{\nu !}\, d_1^{(\nu)}(\nu_n),\quad
d_{\nu n}:=\frac{1}{(\nu+m_n)!}\, d^{(\nu+m_n)}(\nu_n),
\quad \nu=\overline{0,m_n-1}.
$$
In particular, $M_{n+m_n-1}=-d^1_{0n}/d_{0n}.$
If $m_n=1$ (i.e. $n\in I\setminus I'$), then
$$
M_n=-\di\frac{d_1(\nu_n)}{\dot d(\nu_n)},
\qquad \dot d(\rho):=\frac{d}{d\rho}\,d(\rho).                              \eqno(22)
$$
The sequence $\{M_n\}_{n\in\Lambda}$ is called the Weyl sequence. The data
${\cal D}=\{\nu_n, M_n\}_{n\in\Lambda}$ are called the spectral data.
We note that the specification of the spectral data ${\cal D}$ uniquely
determines the Weyl-type function $M(\rho)$ (see [16]).

We consider the following auxiliary inverse problem which is called IP-0.

\smallskip
{\bf IP-0. } Given the spectral data ${\cal D}=\{\nu_n, M_n\}_{n\in\Lambda}$
and $\ga_j,\; j=\overline{1, N-1},$
construct $p(x), q(x), x\in (0, T),$ $\eta'_j, \eta_j,\;j=\overline{1, N-1}.$

\smallskip
Let us prove the uniqueness theorem for IP-0. For this purpose together
with $B$ we consider a boundary value problem $\tilde B$ of the same form
but with different potentials $\tilde p(x), \tilde q(x),$  and different
coefficients of the boundary and the jump conditions. We agree that if
a sertain symbol $\theta$ denotes an object related to $B,$ then
$\tilde\theta$ will denote the analogous object related to $\tilde B.$

\smallskip
{\bf Theorem 1. }{\it If ${\cal D}=\tilde{\cal D},\; \ga_j=\tilde\ga_j,\;
j=\overline{1,N-1},$ then $p=\tilde p, q=\tilde q,$ $\eta'_j=\tilde\eta'_j,
\eta_j=\tilde\eta_j$, $j=\overline{1,N-1}.$}

\smallskip
{\it Proof. } Consider the functions
$$
\left.\begin{array}{c}
P_{11}(x,\rho)=\Phi(x,\rho)\tilde S'(x,\rho)-S(x,\rho)\tilde\Phi'(x,\rho),\\[3mm]
P_{12}(x,\rho)=S(x,\rho)\tilde\Phi(x,\rho)-\Phi(x,\rho)\tilde S(x,\rho).
\end{array}\right\}                                                           \eqno(23)
$$
Using (10), we calculate
$$
\left.\begin{array}{c}
S(x,\rho)=P_{11}(x,\rho)\tilde S(x,\rho)+P_{12}(x,\rho)\tilde S'(x,\rho),\\[3mm]
\Phi(x,\rho)=P_{11}(x,\rho)\tilde\Phi(x,\rho)+P_{12}(x,\rho)\tilde\Phi'(x,\rho).
\end{array}\right\}                                                           \eqno(24)
$$
According to (8) and (9), for each fixed $x,$ the functions
$P_{11}(x,\rho)$ and $P_{12}(x,\rho)$ are meromorphic in $\rho$
with poles at the points $\nu_n$ and $\tilde\nu_n.$ Denote
$G_{\delta}^0=G_{\delta}\cap{\tilde G}_{\delta}.$ By virtue
of (19), (20) and (23) we get
$$
|P_{12}(x,\rho)|\le C|\rho|^{-1},\quad
|P_{11}(x,\rho)|\le C,\quad \rho\in G_{\delta}^0.                            \eqno(25)
$$
It follows from (8) and (23) that
$$
P_{11}(x,\rho)=C(x,\rho)\tilde S'(x,\rho)-S(x,\rho)\tilde C'(x,\rho)
+(M(\rho)-\tilde M(\rho))S(x,\rho)\tilde S'(x,\rho),
$$
$$
P_{12}(x,\rho)=S(x,\rho)\tilde C(x,\rho)-C(x,\rho)\tilde S(x,\rho)
-(M(\rho)-\tilde M(\rho))S(x,\rho)\tilde S(x,\rho).
$$
Since ${\cal D}=\tilde{\cal D},$ it follows from (21) that
for each fixed $x,$ the functions $P_{1k}(x,\rho)$ are entire in $\rho.$
Together with (25) this yields $P_{12}(x,\rho)\equiv 0,\;
P_{11}(x,\rho)\equiv A(x),$ where the function $A(x)$ does not depend
on $\rho.$ Using (24) we derive
$$
S(x,\rho)\equiv A(x)\tilde S(x,\rho),\;\Phi(x,\rho)\equiv A(x)\tilde\Phi(x,\rho). \eqno(26)
$$
Taking (15)-(17) and (23) into account, we obtain
for each fixed $x\in(b_j, b_{j+1})$:
$$
P_{11}(x,\rho)=\frac{1}{2}
\Big(\frac{\xi_1^{\pm},\ldots,\xi_j^{\pm}}{\tilde\xi_1^{\pm},\ldots,\tilde\xi_j^{\pm}}
+\frac{\tilde\xi_1^{\pm},\ldots,\tilde\xi_j^{\pm}}{\xi_1^{\pm},\ldots,\xi_j^{\pm}}\Big)
[1],\quad \rho\in\Pi_{\de}^{\pm},\; |\rho|\to\iy.
$$
Therefore, the function $A(x)$ is piecewise constant (a step-function).
Together with (26) this yields $q(x)=\tilde q(x), p(x)=\tilde p(x),$ $x\in(0,T).$
It follows from (26) that
$$
\frac{\tilde S(x,\rho)}{S(x,\rho)}=\frac{\tilde\Phi(x,\rho)}{\Phi(x,\rho)},
$$
and consequently, $(\xi_j^{\pm})^2=(\tilde\xi_j^{\pm})^2$, $j=\overline{1, N-1}.$
Then $\eta'_j=\tilde\eta'_j$, $j=\overline{1, N-1},$ and consequently,
$\xi_j^{\pm}=\tilde\xi_j^{\pm}$, $j=\overline{1, N-1},$
$S(x,\rho)\equiv \tilde S(x,\rho),$ $\Phi(x,\rho)\equiv \tilde\Phi(x,\rho),$
$\eta_j=\tilde\eta_j$, $j=\overline{1, N-1}.$ Theorem 1 is proved.

Using the method of spectral mappings [3], one can obtain a constructive
procedure for the solution of IP-0 (see [3] for details).\\

{\bf 4. Solition of the Inverse problem 1.}\\

Let $a(\rho), d(\rho)$ and $\Omega$ be given. Note that $\al, \be, \ga_j,
j=\overline{1,N-1}$ are known a priori.

The solution of the Inverse problem 1 is constructed as follows.

First we calculate the zeros ${\cal V}=\{\nu_n\}_{n\in\Lambda}$ of $d(\rho).$
Using the asymptotics (18) we find $(\xi_1^{\pm}\ldots \xi_{N-1}^{\pm})$ and $h'$.
Taking (5) into account, we construct $D(\rho).$ According to (6) we calculate
$Q^2(\nu_n)=D^2(\nu_n)-4\al\be,$ and
$$
Q(\nu_n)=\om_n\sqrt{D^2(\nu_n)-4\al\be},\quad n\in I.                         \eqno(27)
$$
Here and below we agree that if $z=|z|e^{i\xi},\; \xi\in[0,2\pi),$ then
$\sqrt{z}=|z|^{1/2}e^{i\xi/2}.$ We construct $\om_{n0}$ by
$$
\om_{n0}=\frac{1}{2\al}\Big(D(\nu_n)+Q(\nu_n)\Big),                           \eqno(28)
$$
since $\om_{n0}=d_1(\nu_n)=C(T,\nu_n)=\vv(T,\nu_n).$

We construct the Weyl sequence $\{M_n\}_{n\in\Lambda}$ as follows:\\
{\it Case 1. } Let $n\in I\setminus I'$ (i.e. $m_n=1$). Then, in view of (22),
$$
M_n=-\frac{\omega_{n0}}{\dot d(\nu_n)}.                                       \eqno(29)
$$
{\it Case 2. } Let $n\in I_1$ (i.e. $m_n>1,\; \omega_n\ne 0$).
Then it follows from (7) that
$$
(\dot Q(\rho)Q(\rho))^{(\nu-1)}_{|\rho=\nu_n}=
(\dot D(\rho)D(\rho))^{(\nu-1)}_{|\rho=\nu_n},\quad \nu=\overline{1,m_n-1}.   \eqno(30)
$$
Using (30) we find $Q^{(\nu)}(\nu_n),\; \nu=\overline{1,m_n-1}.$ Since
$$
d_1(\rho)=\frac{1}{2\al}\Big(D(\rho)+Q(\rho)\Big)-(i\rho h'+h)d(\rho),
$$
we construct $\omega_{n\nu},\; \nu=\overline{1,m_n-1}$ by the formula
$$
\omega_{n\nu}=\frac{1}{2\al}
\Big(D^{(\nu)}(\nu_n)+Q^{(\nu)}(\nu_n)\Big),\quad \nu=\overline{1,m_n-1}.    \eqno(31)
$$
{\it Case 3.} Let $n\in I_0$ (i.e. $m_n>1,\; \omega_n=0$).
Then $\omega_{n\nu},\; \nu=\overline{1,m_n-1}$ are given a priori.

\smallskip
Thus, we have constructed the Weyl sequence $\{M_n\}_{n\in\Lambda}.$
By solving the auxiliary inverse problem IP-0 we find
$p(x), q(x), x\in (0, T),$ $\eta'_j, \eta_j,\; j=\overline{1, N-1},$
and then we calculate the coefficient $h,$ using (4).
Thus, the following theorem holds.

\smallskip
{\bf Theorem 2. }{\it The specification of $a(\rho), d(\rho)$ and $\Omega$
uniquely determines the boundary value problem $B.$ The solution of Inverse
problem 1 can be found by the following algorithm.}

\smallskip
{\bf Algorithm. }{\it Given $a(\rho), d(\rho)$ and $\Omega$.\\
1) Calculate zeros ${\cal V}=\{\nu_n\}_{n\in\Lambda}$ of $d(\rho).$\\
2) Find $h'$ using (18).\\
3) Construct $D(\rho)$ by (5).\\
4) Calculate $Q^2(\nu_n)=D^2(\nu_n)-4\al\be,$ and $Q(\nu_n)$ by (27).\\
5) Construct $\om_{n0}$ by (28). \\
6) Construct $Q^{(\nu)}(\nu_n),\; n\in I_1,\; \nu=\overline{1,m_n-1}$ using (30).\\
7) Find $\om_{n\nu},\;n\in I_1,\; \nu=\overline{1,m_n-1},$ via (31).\\
8) Calculate the Weyl sequence $\{M_n\}_{n\ge 1}$ using (29) and the
recurrent formula
$$
M_{n+m_n-1-\nu}=-\frac{1}{d_{0n}}\Big(d^1_{\nu n}+\sum_{k=0}^{\nu-1}
M_{n+m_n-1-k}d_{\nu-k,n}\Big),\quad n\in I,\; \nu=\overline{0,m_n-1},\;
d^1_{\nu n}:=\frac{1}{\nu !}\,\om_{n\nu}.
$$
9) Find $p(x), q(x), x\in(0,T)$ and $\eta'_j, \eta_j,\; j=\overline{1, N-1},$
by solving the inverse problem IP-0.\\
10) Calculate the coefficient $h$.}

\smallskip
Similarly, one can solve the following inverse problem.

\smallskip
{\bf Inverse problem 2. }
Given ${\cal P}, {\cal V}$ and $\Omega,$ construct $B.$

\bigskip
{\bf Acknowledgment.} This work was supported by Grant 1.1436.2014K of the Russian
Ministry of Education and Science and by Grant 13-01-00134 of Russian Foundation for
Basic Research.

\begin{center}
{\large REFERENCES}
\end{center}

\begin{enumerate}
\item[{[1]}] Marchenko V.A., Sturm-Liouville operators and their applications.
     "Naukova Dumka",  Kiev, 1977;  English  transl., Birkh\"auser, 1986.
\item[{[2]}] Freiling G. and Yurko V.A., Inverse Sturm-Liouville Problems
     and their Applications. NOVA Science Publishers, New York, 2001.
\item[{[3]}] Yurko V.A. Method of Spectral Mappings in the Inverse Problem
     Theory, Inverse and Ill-posed Problems Series. VSP, Utrecht, 2002.
\item[{[4]}] Marchenko V.A. and Ostrovskii I.V., A characterization of the
     spectrum of the Hill operator, Mat.Sb. 97 (1975), 540-606;
     English transl., Math.USSR-Sb. 26 (1975), 4, 493-554.
\item[{[5]}] Yurko V.A., An inverse problem for second order differential
     operators with regular boundary conditions, Mat. Zametki 18, no. 4 (1975),
     569-576; English transl. in Math. Notes 18 (1975), 3-4, 928-932.
\item[{[6]}] Yurko V.A., The inverse spectral problem for differential operators
     with nonseparated boundary  conditions. J. Math. Anal. Appl. 250 (2000), 266-289.
\item[{[7]}] Gasymov M., Guseinov I.M. and Nabiev I.M., An inverse problem for the
     Sturm-Liouville operator with nonseparated self-adjoint boundary conditions.
     Sib. Mat. Zh. 31, no. 6 (1990), 46-54; English transl. in Siberiam Math. J.
     31, no.6 (1990), 910-918.
\item[{[8]}] Collatz, L. Eigenwertaufgaben mit technischen Anwendungen. Akad.
     Verlagsgesellschaft Geest \& Portig, Leipzig, 1963.
\item[{[9]}] Mennicken R. and M\"oller M. Non-self-adjoint boundary value
     problems. North-Holland Mathematic Studies, vol. 192, Amsterdam,
     North-Holland, 2003.
\item[{[10]}] Shkalikov A.A. Boundary problems for ordinary problems for
     differential equations with parameter in the boundary conditions.
     J. Sov. Math. 33 (1986), 1311-1342; translation from Tr. Semin. im.
     I.G. Petrovskogo 9 (1983), 190-229.
\item[{[11]}] Tretter Ch. Boundary eigenvalue problems with differential
     equations $N\eta=\la P\eta$ with $\la$-- polynomial boundary conditions.
     J. Differ. Equ. 170 (2001), 408-471.
\item[{[12]}] Gasymov M.G. and Gusejnov G.S. Determination of diffusion operators
     from the spectral data. DAN Azer. SSR 37, no. 2, (1981) 19-23.
\item[{[13]}] Yurko V.A. An inverse problem for pencils of differential
     operators. Matem. Sbornik 191 (2000), no. 10, 137-160 (Russian);
     English transl. in Sbornik: Mathematics 191 (2000), no. 10, 1561-1586.
\item[{[14]}] Guseinov I. and Nabiev I. The inverse spectral problem for pencils
     of differential operators. Sb. Math. 198, no.11 (2007), 1579-1598;
     transl. from Mat.Sb. 198, no.11 (2007), 47-66.
\item[{[15]}] Buterin S.A. and Yurko V.A. An inverse spectral problem for pencils
     of differential operators on a finite interval. Vestnik Bashkir. Uni.
     no.4 (2006), 1-7 (in Russian).
\item[{[16]}] Yurko V.A. Inverse problems for non-selfadjoint quasi-periodic differential
     pencils. Analysis and Mathematical Physics, vol.2, no.3 (2012), 215-230.
\item[{[17]}] Yurko V.A. On an inverse spectral problem for differential
     operators on a hedgehog-type graph. Doklady Akad. Nauk 425, no.4 (2009),
     466-470; English transl: Doklady Mathematics 79, no.2 (2009), 250-254.
\item[{[18]}] Yurko V.A. Inverse problems for Sturm-Liouville operators
     on bush-type graphs. Inverse Problems 25, no.10 (2009), 105008, 14pp.
\item[{[19]}]  Yurko V.A., Choque Rivero A. and Karlovich Yu. An inverse problem
     for differential operators on hedgehog-type graphs with general matching conditions. Communications in Mathematical Analysis 17, no.2 (2014), 98-107.
\end{enumerate}

\vspace{0.2cm}

\begin{tabular}{ll}
Name:             &   Yurko, Vjacheslav  \\
Place of work:    &   Department of Mathematics, Saratov State University \\
{}                &   Astrakhanskaya 83, Saratov 410012, Russia \\
E-mail:           &   yurkova@info.sgu.ru\\
\end{tabular}

\end{document}